\documentclass[a4paper,fleqn]{cas-dc}
\usepackage[numbers]{natbib}

\usepackage{graphicx,amssymb,amsmath}
\usepackage{rotating}

\newtheorem{algorithm}{Algorithm}
\newtheorem{theorem}{Theorem}
 \newtheorem{lemma}{Lemma}
 \newdefinition{rmk}{Remark}
 \newproof{pf}{Proof}
 \newproof{pot}{Proof of Theorem \ref{fla-thm}}
 
\def\rit{\mathbb{R}}
\def\endprf{\qed}
\def\personalitem{{\boldmath$\cdot$\;}}
\def\fmodel{\check{f}}
\def\bundleset{\mathcal B}

\def\flow{f_{\text{low}}}
\def\fbest{f_{\text{best}}}

\begin{document}
\shorttitle{Fast proximal algorithms}
\title [mode = title]{Fast proximal algorithms for nonsmooth convex optimization}

\author[]{Adam Ouorou}

\address{
Orange Labs Research,
44 avenue de la R\'epublique, 
92300 Chatillon, France.}
\ead{adam.ouorou@orange.com}

\begin{abstract}
In the lines of our approach in \cite{Ouorou2019}, where we exploit Nesterov fast gradient concept \cite{Nesterov1983} to the Moreau-Yosida
  regularization of a convex function, we devise  new proximal algorithms for nonsmooth convex optimization. 
  These algorithms need no bundling mechanism 
  to update the stability center while preserving
  the complexity estimates established in~\cite{Ouorou2019}. We report some preliminary computational results on some academic test problem to
  give a first estimate of their performance in relation with the classical proximal bundle algorithm.
\end{abstract}


\maketitle

\section{Introduction}
We are interested in minimizing a nonsmooth convex function $f:\;\rit^n\to\rit$, over a  nonempty convex compact subset $S$ of $\rit^n$.
We denote by $f^*$ the optimal value of this problem and $x^*$ an optimal solution. 
Having generated a number of test points 
$y^i\in S,\;i=1,\ldots$ 
with the corresponding function values $f(y^i)$ and subgradients $g^i\in\partial f(y^i)$ via an \emph{oracle} (for $f$) to form the bundle $\bundleset:\{(y^i,f(y^i), g^i)\}$,  the function 
\begin{equation}\label{fmodelk}
\fmodel_{\bundleset}(x)  = \max\{f(y^i)+\langle g^i,\; 
x-y^i\rangle,\;i\in\bundleset\},
\end{equation}
is a piecewise cutting-plane model for $f$, which underestimates $f$,
i.e. for any $x\in S$, $\fmodel_{\bundleset}(x)\le f(x)$. We use the shortcut $i\in\bundleset$ for $(y^i,f(y^i),g^i)\in\bundleset$. 
Let  $F_\mu$ be the Moreau-Yosida regularization of $f$ w.r.t. some $\mu>0$ assumed fixed in the sequel. The function $F_\mu$ is given by
$$
F_\mu(x)=\min\limits_{z\in S}\left\{f(z)+\dfrac{\mu}{2}\|z-x\|^2\right\}.
$$
Minimizing $f$ is equivalent to minimizing $F_\mu$.
Exploiting the fact that $F_\mu$ is convex and differentiable, it is proposed in \cite{Ouorou2019}, 
to apply the concept of fast gradient method \cite{Nesterov1983,Nesterov2004} to $F_\mu$ for the  minimization of $f$.
This results in the following scheme, starting from any $x^0=y^0\in\rit^n$,
\begin{equation}\label{nagd-ipa}
  \begin{array}{l}
      y^{k+1}=\arg\min\limits_{x\in S}\left\{f(x)+\dfrac{\mu}{2}\|x-x^k\|^2\right\}=x^k-\dfrac{1}{\mu}\nabla F_\mu(x^k),\\
      \\[-2mm]
      x^{k+1}=y^{k+1}+\alpha_k(y^{k+1}-y^k),\;\alpha_k=\lambda_{k+1}^{-1}(\lambda_k-1),
      \end{array}
\end{equation}
where $\{\lambda_k\}$ is the Nesterov' sequence  defined by
$$
\lambda_0=1,\;\;\lambda_{k+1}=\dfrac{1+\sqrt{1+4\lambda_{k}^2}}{2},\;k\ge 0.
$$
This sequence has the following properties
\begin{equation}\label{l-prop}
\lambda_{k-1}^2=\lambda_k^2-\lambda_k,\;k\ge 1,\;\lambda_k^2=\sum\limits_{i=0}^k\lambda_i,\;\lambda_k\ge\frac{k+2}{2},\;k\ge 0.
\end{equation}
The above scheme generates a sequence $\{y^k\}$ of approximations to an optimal point of the considered problem, and a second sequence
$\{x^k\}$ of \emph{stability centers}, different from the former. 
It is possible to use another update for $x^{k+1}$, as proposed by G\"uler in~\cite{Guler1992}, 
\begin{equation}\label{guler}
x^{k+1}  =y^{k+1}+\alpha_k(y^{k+1}-y^k)+\beta_k(y^{k+1}-x^k),
\end{equation}
where
\begin{equation}\label{beta}
\beta_k=\lambda_k\lambda_{k+1}^{-1}.
\end{equation}

Computing exactly $y^{k+1}$, the proximal point of the stability center $x^k$ is out of reach in practice. In \cite{Ouorou2019}, we compute 
an approximate solution through a sequence of quadratic subproblems 
$$
z^j=\arg\min\limits_{x\in S}\left\{\fmodel_{\bundleset_j}(x)+\dfrac{\mu}{2}\|x-x^k\|^2\right\},\;\;j=1,\ldots,
$$
As the bundle $\bundleset_j$ grows, $\{z^j\}_{j\ge 0}$ tends to $y^{k+1}$. An approximate proximal point of $x^k$ 
is identified when the condition 
\begin{equation}\label{bundle-stop}
f(z^j)-\fmodel_{\bundleset_j}(z^j)\le\varepsilon_k,
\end{equation}
is satisfied for some positive tolerance $\varepsilon_k$, in which case  
$y^{k+1}$ is set to $z^j$ (we keep the same notation as for  the exact proximal point).
Then, the next stability center $x^{k+1}$ is updated using this approximation 
in place of the exact proximal point 
in \eqref{nagd-ipa} or \eqref{guler}. There are two versions of the above outlined algorithm FPBA,
that we denote by FPBA1 and FPBA2 (FPBA stands for Fast Proximal Bundle Algorithm). FPBA1 uses the rule in \eqref{nagd-ipa} to update the next prox-center $x^{k+1}$
 while FPBA2 uses the \emph{momentum term} proposed by G\"uler with the sequence $\{\beta_k\}$, cf \eqref{guler}. 
The complexity estimates of the two algorithms are given respectively as
\begin{equation}\label{nesterov-estimate}
f(y^k)-f^*\le\frac{2\mu \|x^0-x^*\|^2}{(k+1)^2}+\vartheta_k,\;\;k\ge 1,
\end{equation}
for FPBA1, and
\begin{equation}\label{guler-estimate}
f(y^k)-f^*\le\frac{\mu \|x^0-x^*\|^2}{(k+1)^2}+\vartheta_k,
\;\;k\ge 1,
\end{equation}
for FPBA2, where $x^*$ is any optimal solution and
\begin{equation}\label{acc-errors}
  \vartheta_k=\lambda_{k-1}^{-2}\sum\limits_{i=0}^{k-1}\lambda_{i}^2\varepsilon_i
\end{equation}
is the accumulation of errors at step~$k$, see Theorems~3.1 and~3.2 in \cite{Ouorou2019}.

In this paper, we first propose a variant of the algorithm FPBA in \cite{Ouorou2019} that uses no bundling mechanism to update the stability center $x^k$ 
and second, take inspiration from \cite{LNN95,OlSo2016} and the recent survey on bundle methods \cite{Frangioni2020}
we propose two other proximal algorithms in the lines of our 
approach. For ease of exposition, we limit our development to the case with $\beta_k=0,\;k\ge 0$. 
The results for the other case 
can be obtained  conjointlty with the development and arguments used in \cite{Ouorou2019}.
We assume that $S$ is simple enough to allow solving easily all the linear and quadratic subproblems in the paper.
\section{Fast proximal cutting plane algorithm} 
In  Section~4 of \cite{Ouorou2019}, an analysis of the accumulated errors \eqref{acc-errors} 
shows that
one can tolerate large errors in early iterations but require smaller and smaller errors in the progress of the algorithms. Based on this, on may be content with  
only one quadratic subproblem at each step~$k$, obtaining the proximal point of $x^k$ w.r.t to $\fmodel_{\bundleset_k}$, 
which is an approximation of the exact proximal point of $x^k$ w.r.t. $f$ with error $\varepsilon_k=f(y^{k+1})-\fmodel_{\bundleset_j}(y^{k+1})\ge 0$. 
There is no need to distinguish between serious and null steps as in FPBA or classical proximal bundle algorithms. 
The resulting variant of FPBA, which we term as Fast Proximal Cutting Plane Algorithm (FPCPA), is as follows.
\begin{algorithm}\label{newfpba}
\begin{enumerate}
\item[] 
\item[0.]\label{fpba-stepi}  Choose $x^0=y^0\in \rit^n$ 
  and the sequence $\{\beta_k\}_{k\ge 0}$.  Set $k=0$.
\item\label{fpba-step0} Compute $f(y^k),\;
  g^k\in\partial f(y^k)$ and update $\mathcal B_k$.
\item\label{fpba-step1}  If $g^{k}=0$, terminate.
\item\label{fpba-step2} Compute
\begin{equation}\label{bundle-qp}
  y^{k+1}=\arg\min\limits_{x\in S}\left\{\fmodel_{\bundleset_k}(x)+\dfrac{\mu}{2}\|x-x^k\|^2\right\},
\end{equation}
and $x^{k+1}  =y^{k+1}+\alpha_k(y^{k+1}-y^k)+\beta_k(y^{k+1}-x^k)$.
\item Set $k=k+1$ and go to Step \ref{fpba-step0}.
\end{enumerate}
\end{algorithm}
In this algorithm, the choice of  
the sequence $\{\beta_k\}$ as $\beta_k=0,\;k\ge 0$ or $\beta_k=\lambda_k\lambda_{k+1}^{-1},\;k\ge 0$, 
results in two versions of the algorithm, which we denote respectively by FPCPA1 and FPCPA2. 
They preserve respectively the complexity estimates \eqref{nesterov-estimate} and \eqref{guler-estimate}. 
It is possible to use a proximity parameter that depends on $k$ with the same complexity estimates,
provided that $\mu_0=\mu$ and $\mu_k\le\mu_{k-1},\;k\ge 1$,
see Proposition~3.1 in \cite{Ouorou2019}. With a dynamic setting of the proximity parameter, Algorithm~\ref{newfpba} appears as an
implementable version of the inertial proximal algorithm \cite{AtCa2018}. The convergence of this algorithm may be derived from that of
Algorithm~\ref{fast-dbl}  below.

\section{Fast level algorithm}
Define the \emph{level} $l_k$ by 
\begin{equation}\label{level-def}
l_k=\kappa \flow^k+(1-\kappa)\fbest^k =\fbest^k-\kappa\Delta_k, 
\end{equation}
where
\begin{itemize}
\item[\personalitem] $0<\kappa< 1$ is the \emph{level parameter},
\item[\personalitem] $\fbest^k$ is the best objective value found at step~$k$, 
\item[\personalitem] $\flow^k$ is a finite lower bound on $f^*$, 
\item[\personalitem] $\Delta_k=\fbest^k-\flow^k\ge 0$. 
\end{itemize}
By interpreting the term $\fmodel_{\bundleset_k}(x)$ in \eqref{bundle-qp} as the dualization of a constraint $\fmodel_{\bundleset_k}(x)\le l_k$,
an alternative to \eqref{bundle-qp}, consists in projecting $x^k$ on the $l_k$-level set of $\fmodel_{\bundleset_k}$, see \cite{Frangioni2020,HUL1993}. 
The corresponding algorithm, denoted by FLA (for Fast Level Algorithm) is as follows.
\begin{algorithm}\label{fla}
\begin{enumerate}
\item[] 
\item[0.]\label{fla-stepi}  Choose $x^0=y^0\in \rit^n$ and the sequence $\{\beta_k\}_{k\ge 0}$. Set $k=0$.
\item\label{fla-step1} Compute $f(y^{k})$ and $g^{k}\in\partial f(y^{k})$ and update $\mathcal B_k$.
\item\label{fla-step2} Update $\fbest^k,\;\flow^k$. Set $\Delta_k=\fbest^k-\flow^k$ and $l_k=\fbest^k-\kappa\Delta_k$.
\item\label{fla-step3} If $\Delta_k\le\varepsilon$ or $g^{k}=0$, stop.
\item\label{fla-step4}  Compute 
\begin{equation}\label{fla-qp}
y^{k+1}=\arg\min\limits_{x\in S}\left\{\frac{1}{2}\|x-x^k\|^2:\;\fmodel_{\bundleset_k}(x)\le l_k\right\}
\end{equation}
and $x^{k+1}  =y^{k+1}+\alpha_k(y^{k+1}-y^k)+\beta_k(y^{k+1}-x^k)$.
\item Set $k=k+1$ and loop to Step~\ref{fla-step1}.
\end{enumerate}
\end{algorithm}
The convergence property of this algorithm is given below.
\begin{theorem}\label{fla-thm}
  For the sequence $\{y^k\}$ generated by Algorithm~\ref{fla} with $\beta_k=0,\;k\ge 0$, we have
    $$
    f(y^k)-f^*\le\frac{2\mu\|x^0-x^*\|^2}{t_0(k+1)^2}+\vartheta_k,\;k\ge 1,
    $$
    where $t_k$ is the optimal dual solution of \eqref{fla-qp-eq}.
  \end{theorem}
  The proof is given after that of the next algorithm. 
\section{Fast doubly stabilized algorithm}
In this section, taking inspiration from \cite{OlSo2016}, we propose an algorithm with the aim to leverage 
the good features of the two previous ones by combining the quadratic problems \eqref{bundle-qp} and 
\eqref{fla-qp} into a single quadratic subproblem as follows
$$
\min\limits_{x\in S}\left\{\fmodel_{\bundleset_k}(x)+\dfrac{\mu_k}{2}\|x-x^k\|^2:\; \fmodel_{\bundleset_k}(x)\le l_k\right\},
$$
or equivalently
\begin{equation}\label{dbl-qp}
  \min\limits_{(x,r)\in S\times\rit}\left\{r+\frac{\mu_k}{2}\|x-x^k\|^2: \; \fmodel_{\bundleset_k}(x)\le r,\;r\le l_k\right\}.
\end{equation}
For a reason to be apparent shortly, here the proximity parameter  needs to depend on $k$. The resulting
algorithm is as follows, we term it as Fast Doubly Stabilized Algorithm (FDSA for short), keeping the wording of \cite{OlSo2016}.
\begin{algorithm}\label{fast-dbl}
\begin{enumerate}
\item[] 
\item[0.]\label{fast-dbl-stepi}  Choose $x^0=y^0\in \rit^n$ and the sequence $\{\beta_k\}_{k\ge 0}$. Set $k=0$.
\item\label{fast-dbl-step1} Compute $f(y^{k})$ and $g^{k}\in\partial f(y^{k})$ and update $\mathcal B_k$.
\item\label{fast-dbl-step2} Update $\fbest^k,\;\flow^k$. Set $\Delta_k=\fbest^k-\flow^k$ and $l_k=\fbest^k-\kappa\Delta_k$.
\item\label{fast-dbl-step3} If $\Delta_k\le\varepsilon$ or $g^{k}=0$, stop.
\item\label{fast-dbl-step4}  Compute the $x$-solution $y^{k+1}$ of \eqref{dbl-qp} and set
 $x^{k+1}  =y^{k+1}+\alpha_k(y^{k+1}-y^k)+\beta_k(y^{k+1}-x^k)$. 
\item Set $k=k+1$ and loop to Step~\ref{fast-dbl-step1}.
\end{enumerate}
\end{algorithm}
Its convergence is given by the next result.  
\begin{theorem}\label{dbl-thm}
  Given some $\mu>0$, assume that the sequence $\{\mu_k\}$ satisfies $\mu_0=\mu$ and $\mu_kt_{k-1}\le\mu_{k-1}t_k$ for $k\ge 1$. Then, for the sequence $\{y^k\}$
generated by Algorithm~\ref{fast-dbl} with $\beta_k=0,\;k\ge 0$, we have 
$$
f(y^k)-f^*\le\frac{2\mu\|x^0-x^*\|^2}{t_0(k+1)^2}+\vartheta_k,\;k\ge 1,
$$
where $t_k$ is the optimal dual solution associated with the constraint $\fmodel_{\bundleset_k}(x)\le r$ in \eqref{dbl-qp}
and $\vartheta_k$ is given by \eqref{acc-errors}.
\end{theorem}
\begin{pf}
The proof uses the arguments of Lemma 3.1, 3.2 and Theorem 3.1. For the paper to be self-contained, we povide a complete proof.
  
The KKT conditions for \eqref{dbl-qp} imply that there exist
$p_f^k\in\partial\fmodel_{\bundleset_k}(y^{k+1}),\;p_S^k\in\partial\mathcal I_S(y^{k+1})$ and real numbers $t_k,\;\tau_k\ge 0$ such that 
\begin{equation}\label{dbl-kkt}
\begin{array}{l}
 \mu_k( y^{k+1}-x^k)+t_k(p_f^k+p_S^k)=0,\;t_k\lbrack\fmodel_{\bundleset_k}(y^{k+1})-r^k\rbrack=0,\\
  \\[-2mm]
  1-t_k+\tau_k=0,\;\tau_k(r^k-l_k)=0,
  \end{array}
\end{equation}
where $r^k$ is the $r$-solution of \eqref{dbl-qp}. 
These conditions imply that $t_k=\tau_k+1\ge 1$ and 
\begin{equation}\label{dbl-grad}
  \gamma_k(x^k-y^{k+1})=p_f^k+p_S^k\;\text{where}\;\gamma_k=t_k^{-1}\mu_k\;(\le\mu_k).
\end{equation}
Recall that $\partial I_S(x)$ is the normal cone of $S$ at $x$ i.e.
$$
\partial I_S(x)=\{y\in\rit^n:\;\langle y,z-x\rangle\le 0,\;z\in S\}.
$$
We have for any $x\in S$,
$$
\begin{array}{lcl}
 \langle p_f^k+p_S^k,x-y^{k+1}\rangle
  &=&\langle p_f^k,x-y^{k+1}\rangle+\langle p_S^k,x-y^{k+1}\rangle\\
  \\[-2mm]
  & \le & \langle p_f^k,x-y^{k+1}\rangle.
  \end{array}
  $$
  Therefore, as $p_f^k\in\partial\fmodel_{\bundleset_k}(y^{k+1})$, we get for any $x\in S$,
  $$
  \fmodel_{\bundleset_k}(y^{k+1})+\langle p_f^k+p_S^k,x-y^{k+1}\rangle\le \fmodel_{\bundleset_k}(x)\le f(x).
  $$
  and
 $f(y^{k+1})+\langle p_f^k+p_S^k,x-y^{k+1}\rangle-\varepsilon_k\le f(x)$, 
  with $\varepsilon_k=f(y^{k+1})-\fmodel_{\bundleset_k}(y^{k+1})$. 
  In other words, $p_f^k+p_S^k\in\partial_{\varepsilon_k} f(y^{k+1})$.
Using \eqref{dbl-grad}, we have for any $x\in S$,
 \begin{equation}\label{dbl-sub-ineq}
f(x)\ge f(y^{k+1})+\gamma_k\langle x^k-y^{k+1},x-y^{k+1}\rangle-\varepsilon_k.
\end{equation}
Let $\delta_k=f(y^k)-f^*$. Taking $x=y^k(\in S)$ in \eqref{dbl-sub-ineq} and multiplying the resulting inequality with 
$\lambda_{k}-1$ give
$$
(\lambda_{k}-1)(\delta_k-\delta_{k+1})\ge \gamma_k\langle x^k-y^{k+1},\lambda_k(y^k-y^{k+1})+y^{k+1}-y^k\rangle-\lambda_k\varepsilon_k.
$$
We add this inequality with the one resulting from \eqref{dbl-sub-ineq} with  $x=x^*(\in S)$ and get
$$
(\lambda_{k}-1)\delta_k-\lambda_k\delta_{k+1}\ge \gamma_k\langle x^k-y^{k+1},\lambda_k(y^k-y^{k+1})+x^*-y^k\rangle-\lambda_k\varepsilon_k
$$
Now, multiplying the above inequality by $\lambda_k$ and using the first relation in \eqref{l-prop} yield 
\begin{equation}\label{dbl-delta}
\lambda_{k-1}^2\delta_k-\lambda_k^2\delta_{k+1}\ge \gamma_k\langle u^k,v^k\rangle-\lambda_k^2\varepsilon_k,
\end{equation}  
where
$u^k=\lambda_k(y^{k+1}-x^k)$ and $v^k=\lambda_k(y^{k+1}-y^k)+y^k-x^*$.
For any $u,v\in\rit^n$, we have (parallelogram law)
\begin{equation}\label{p-law}
\langle u,v\rangle  = \frac{1}{2}(\|u\|^2+\|v\|^2-\|u-v\|^2)\ge \frac{1}{2}(\|v\|^2-\|u-v\|^2).
\end{equation}  
Hence
$$
\lambda_{k-1}^2\delta_k-\lambda_k^2\delta_{k+1}\ge \frac{\gamma_k}{2}(\|v^k\|^2-\|v^k-u^k\|^2)-\lambda_k^2\varepsilon_k.
$$
Let $w^k=v^k-u^k=\lambda_k(x^k-y^k)+y^k-x^*,\;k\ge 0$. Then,
$$
\begin{array}{lcl}
  w^{k+1} & = & \lambda_{k+1}(x^{k+1}-y^{k+1})+y^{k+1}-x^*\\
          &\stackrel{\eqref{nagd-ipa}}{=} & (\lambda_k-1)(y^{k+1}-y^k)+y^{k+1}-x^*\\
          &= & v^k,
\end{array}
$$
and
$$
\lambda_{k-1}^2\delta_k-\lambda_k^2\delta_{k+1}\ge \dfrac{\gamma_k}{2}\|w^{k+1}\|^2-\dfrac{\gamma_k}{2}\|w^{k}\|^2-\lambda_k^2\varepsilon_k,
$$
The assumption $\mu_kt_{k-1}\le\mu_{k-1}t_k$ implies $\gamma_k\le\gamma_{k-1}$ and then
$$
\lambda_{k-1}^2\delta_k-\lambda_k^2\delta_{k+1}\ge\dfrac{\gamma_k}{2}\|w^{k+1}\|^2-\dfrac{\gamma_{k-1}}{2}\|w^{k}\|^2-\lambda_k^2\varepsilon_k.
$$
We now sum these inequalities for $i=1,\ldots, k-1$ to get 
\begin{equation}\label{delta-ub}
\begin{array}{lcl}
\lambda_{k-1}^2\delta_k & \le &
\lambda_0\delta_1+\dfrac{\gamma_0}{2}\|w^1\|^2+\sum\limits_{i=1}^{k-1}\lambda_i^2\varepsilon_i-\dfrac{\gamma_{k-1}}{2}\|w^k\|^2\\
  & \stackrel{\lambda_0=1}{\le} &
\delta_1+\dfrac{\gamma_0}{2}\|w^1\|^2+\sum\limits_{i=1}^{k-1}\lambda_i^2\varepsilon_i.
\end{array}
\end{equation}
Using \eqref{dbl-sub-ineq} with $x=x^*$ and $k=0$, we get
$$
\begin{array}{lcl}
  \delta_1 &\le & -\gamma_0\langle x^0-y^{1},\;x^*-y^1\rangle+\varepsilon_0\\
\\[-2mm]
& \stackrel{\eqref{p-law}}{=} & -\dfrac{\gamma_0}{2}\left[\|x^0-y^{1}\|^2+\|y^1-x^*\|-\|x^0-x^*\|^2\right]+\varepsilon_0\\
\\[-2mm]
& \le & -\dfrac{\gamma_0}{2}\|y^1-x^*\|^2+\dfrac{\gamma_0}{2}\|x^0-x^*\|^2
      +\varepsilon_0
\end{array}
$$
Since $w^1=v^0=\lambda_0(y^1-y^0)+y^0-x^*\stackrel{\lambda_0=1}{=}y^1-x^*$, we have
$$
\delta_1+\dfrac{\gamma_0}{2}\|w^1\|^2\le \dfrac{\gamma_0}{2}\|x^0-x^*\|^2+\varepsilon_0,
$$
and from \eqref{delta-ub},
$$
\delta_k \le \dfrac{\gamma_0}{2\lambda_{k-1}^2}\|x^0-x^*\|^2+\lambda_{k-1}^{-2}\sum\limits_{i=0}^{k-1}\lambda_i^2\varepsilon_i.
$$
It remains to use in the first term of the r.h.s. of this inequality, the fact that $\gamma_0=t_0^{-1}\mu$ and $\lambda_{k-1}\ge (k+1)/2$ from~\eqref{l-prop}.
\endprf
\end{pf}

A few comments are in order.
\begin{enumerate}
\item In the same way as Lemma~1 of \cite{OlSo2016},
  it can be shown that the $x$-solution of \eqref{dbl-qp} is either the one of \eqref{bundle-qp} or that of \eqref{fla-qp}.
  Algorithm~\ref{fast-dbl} makes the choice automatically depending on the value of $t_k$ (in fact this choice depends on the proximity and the level
  parameters $\mu_k$  and $\kappa$ (defining $l_k$) which determine $t_k$).
If $y_p^{k+1}$ and $y_l^{k+1}$ denote the respective optimal
solutions of the quadratic problems \eqref{bundle-qp} and \eqref{fla-qp}, we have
$$
y^{k+1}=\left\{
  \begin{array}{ll}
    y_p^{k+1} & \text{if}\;t_k=1 (\tau_k=0),\\
    \\[-2mm]
    y_l^{k+1} & \text{if}\;t_k>1 (\tau_k>0).
  \end{array}
  \right.
$$
Because $t_k>0$, we have $\fmodel_{\bundleset_k}(y^{k+1})=r^k,\;k\ge 0$, 
while $r^k\le l_k\; \text{if}\;t_k=1\;\text{and}\; r^k=l_k \; \text{if}\;t_k>1$.
\item
Discarding the accumulation of errors, the complexity estimate improves slightly compared to the one of
Algorithm~\ref{newfpba} as $t_0\ge 1$, cf \eqref{nesterov-estimate}.
\item
  We get from \eqref{dbl-sub-ineq} and Cauchy-Schwartz inequality, 
  $$
   f(y^{k+1})\le f(x)+\gamma_k\|x^k-y^{k+1}\|\|x-y^{k+1}\|+\varepsilon_k,
   $$
for any $x\in S$. Therefore, if
$$
\gamma_k\|x^k-y^{k+1}\|\le \varepsilon_1\;\text{ and }\;\varepsilon_k\le\varepsilon_2,
$$ 
for some stopping tolerances $\varepsilon_1,\varepsilon_2>0$, then
$$
f(y^{k+1})\le f(x)+\varepsilon_1\|x-y^{k+1}\|+\varepsilon_2,\;\forall x\in S.
$$
We can then consider $y^{k+1}$  as an approximate optimal solution if $\varepsilon_1$ and $\varepsilon_2$ are small enough.
\item
  We can recover the complexity estimate of Algorithm~\ref{newfpba} from Theorem~\ref{dbl-thm}. Indeed,  by replacing \eqref{level-def} with $l_k=+\infty$, we have 
  $\tau_k=0,\;k\ge 0$ and then, $t_k=1$ and $\gamma_k=\mu_k,\;k\ge 0$. Algorithm~\ref{fast-dbl} then reduces to Algorithm~\ref{newfpba}.
  The complexity estimate given in Theorem~\ref{dbl-thm} becomes \eqref{nesterov-estimate}, the one already given for Algorithm~\ref{newfpba},
  with the assumption that now writes $\mu_0=\mu$ and $\mu_k\le\mu_{k-1}$ for $k\ge 1$.
\item 
  For the assumption in Theorem~\ref{dbl-thm} to hold, the proximity parameter needs to depend on $k$.
  An example of  sequence 
  that satisfies this assumption 
  is given by 
  \begin{equation}\label{fdsa-mu}
  \mu_0=\mu\;\text{ and  }\;\mu_k=\gamma_{k-1}=t_{k-1}^{-1}\mu_{k-1},\;k\ge 1,
  \end{equation}
due to the fact that $t_k\ge 1$. 
This rule maintains ($t_{k-1}=1$) or decreases ($t_{k-1}>1$) the proximity parameter for the next step. The sequence is only decreasing then while it may be useful sometimes
to increase the proximity parameter. 
If it was possible to guess $t_k$, an intuitive choice suggested by the assumption would be $\mu_k=t_k\gamma_{k-1}=t_kt_{k-1}^{-1}\mu_{k-1}$. 
The ratio $t_kt_{k-1}^{-1}$ would reflect 
the change between steps $k-1$ and $k$, maintaining ($t_{k-1}=t_k=1$),
increasing ($t_{k-1}<t_k$) or decreasing ($t_{k-1}>t_k$) the proximity parameter accordingly for the next step. Unfortunately, $t_k$ is obtained only
after fixing $\mu_k$ and solving \eqref{dbl-qp}.
\endprf
\end{enumerate}

In the light of the proof of Theorem~\ref{dbl-thm}, we now give that of Theorem~\ref{fla-thm}.
\begin{pot}
It is clear that the unique solution of \eqref{fla-qp} is also the unique solution of the problem
\begin{equation}\label{fla-qp-eq}
\min\limits_{x\in S}\left\{
  \frac{\mu_k}{2}\|x-x^k\|^2:\;\fmodel_{\bundleset_k}(x)\le l_k\right\},
\end{equation}
for any given $\mu_k>0$; we set $\mu_0=\mu$.  
The KKT conditions for this quadratic problem  imply that there exist
$p_f^k\in\partial\fmodel_{\bundleset_k}(y^{k+1}),\;p_S^k\in\partial\mathcal I_S(y^{k+1})$ and $t_k> 0$ (as $x^k\neq y^{k+1}$) 
such that 
$$
\mu_k(y^{k+1}-x^k)+t_k(p_f^k+p_S^k)=0,\;t_k\lbrack\fmodel_{\bundleset_k}(y^{k+1})-l_k\rbrack=0.
$$
Therefore, for any $x\in S$,
  $$
  \fmodel_{\bundleset_k}(y^{k+1})+\langle p_f^k+p_S^k,x-y^{k+1}\rangle 
  \le f(x),
  $$
 or equivalently,
 \begin{equation}\label{sub-ineq}
f(x)\ge f(y^{k+1})+\frac{\mu_k}{t_k}\langle x^k-y^{k+1},x-y^{k+1}\rangle-\varepsilon_k,
\end{equation}
where $\varepsilon_k=f(y^{k+1})-\fmodel_{\bundleset_k}(y^{k+1})$.
The remaing of the proof is similar to that of Theorem~\ref{dbl-thm}, under the same assumption $\mu_0=\mu,\;\mu_kt_{k-1}\le\mu_{k-1}t_k,\;k\ge 1$.
\endprf
\end{pot}
\begin{rmk}
As given in Theorem~\ref{fla-thm}, we cannot state if the complexity estimate Algorithm~\ref{fla} improves or not over the one of
Algorithm~\ref{newfpba}. Since the sequence $\{y^k\}$ generated by Algorithm~\ref{fla} is the same as if we use \eqref{fla-qp-eq} in place of \eqref{fla-qp} in Step~\ref{fla-step4}, 
we conjecture that the complexity estimate of Algorithm~\ref{fla} to be the same as that of Algorithm~\ref{fast-dbl}. Indeed, with an 
appropriate choice of $\mu_k$ and $\kappa$ to have $\tau_k>0$ for all $k\ge 0$, the sequence $\{y^k\}$ generated by Algorithm~\ref{fast-dbl}
is the same as that obtain from Algorithm~\ref{fla}, and then the same complexity estimate
as given by Theorem~\ref{dbl-thm}. 
\endprf
\end{rmk}
\begin{rmk}
The level parameter $\kappa$ does not appears explicitly in the complexity estimate of Algorithm~\ref{fla} as it is for the level bundle algorithms
in \cite{LNN95}. In fact it is hidden in $t_0$ as it influences the level 
and the dual variables of the level constraints.
\endprf
\end{rmk}

We now give the complexity estimates of Algorithms~\ref{fla} and~\ref{fast-dbl} used with the sequence $\{\beta_k\}$ given by \eqref{beta}. The
proof is analogue to that of Theorem~3.2 in \cite{Ouorou2019} with the same arguments used in  the above proofs of Theorems~\ref{fla-thm} and~\ref{dbl-thm}. The main difference
is a better lower bound obtained on the scalar product $\langle u^k,v^k\rangle$ (cf \eqref{dbl-delta}) thanks to the update of $x^{k+1}$ using a second momentum term
proposed by G\"uler intuitively in \cite{Guler1992}. It is shown in \cite{KimFessler2016,KimFessler2018} by Kim and Fessler  that it corresponds to an optimal choice of
parameters obtained through a relaxed \emph{performance estimation problem}
introduced by Drori and Teboulle to optimize first-order algorithms, see \cite{DroriTeboulle2014}.
\begin{theorem}
  Assume that Algorithms ~\ref{fla} and~\ref{fast-dbl} use the sequence \eqref{beta} under the assumption of Theorem~\ref{dbl-thm} on the sequence $\{\mu_k\}$.
  Then, for the sequence $\{y^k\}$ generated, we have
  $$
  f(y^k)-f^*\le\frac{\mu\|x^0-x^*\|^2}{t_0(k+1)^2}+\vartheta_k,\;k\ge 1,
  $$
  where $t_0$ is the (respective) dual solution associated with the constraint $\fmodel_{\bundleset_k}(x)\le w$ in respectively the quadratic
  subproblems \eqref{fla-qp-eq} and \eqref{dbl-qp} with $\mu_0=\mu$,
  and $\vartheta_k$ is given by \eqref{acc-errors}.
\end{theorem}

For the above complexity estimates 
to be meaningful, it is necessary that the accumulation of errors $\vartheta_k$
to not be divergent with the first terms.
\begin{lemma}
The sequence $\{\vartheta_k\}$ is bounded above.
\end{lemma}
\begin{pf}
Recall that
$$
0\le \vartheta_k=\sum\limits_{i=0}^{k-1}\nu_i^k\varepsilon_i\;\text{where}\;\nu_i^k=\lambda_i^2\lambda_{k-1}^{-2},\;i=0,\ldots,k-1.
$$
We can observe that $0\le\nu_i^k\le 1=\nu_{k-1}^k$ and the former errors vanish with their weights as they tend to $0$ when $k$ grows. 
We have 
$\vartheta_{k+1}-\vartheta_k=\varepsilon_k-\lambda_k^{-1}\vartheta_k$,
see Section~4 in \cite{Ouorou2019}. Therefore, 
$\varepsilon_k\le \lambda_k^{-1}\vartheta_k$ implies $\vartheta_{k+1}\le\vartheta_k$. 
From Proposition~4.3 in \cite{CoLe1993}, as the bundle $\bundleset_k$ grows, $f(y^{k+1})$ and $\fmodel_{\bundleset_k}(y^{k+1})$ get closer to each
other i.e. 
$\varepsilon_k=f(y^{k+1})-\fmodel_{\bundleset_k}(y^{k+1})\to 0$ (this means that the last errors vanish as well
with high $k$).
We cannot have $\varepsilon_k>\lambda_k^{-1}\vartheta_k(\ge 0)$ for an infinite number of $k$ as it results in the contradiction $0>0$. 
Therefore, there exists some $k^*$ such that $\varepsilon_k\le \lambda_k^{-1}\vartheta_k$  for $k\ge k^*$ and then
$\vartheta_{k}\le\vartheta_{k-1}\le\ldots\le\vartheta_{k^*+1}\le\vartheta_{k^*}<\infty$, 
 i.e. the sequence $\{\vartheta_k\}_{k\ge k^*}$ is decreasing.
\endprf
\end{pf}
\begin{rmk}
  We finally observe that in the above development, we may replace $\fmodel_k$ by any other 
  lower model $\underline f_k\le f$ and practical
  in the sense that the corresponding subproblems analogue
  to \eqref{bundle-qp}, \eqref{fla-qp} and \eqref{dbl-qp} are easy to solve. In this case, the error at setp $k$ writes 
  $\varepsilon_k=f(y^{k+1})-\underline f_k(y^{k+1})\ge 0$.
\endprf
\end{rmk}
\section{Numerical experiments}
We conducted some preliminary experiments that aim to provide a first look on the
performances of the proposed algorithms as compared with the classical proximal bundle algorithm (CPBA).
The test problems are the one considered in \cite{Ouorou2019} and described in \cite{LukV2000} and the algorithms are implemented using Python 3.5 and Cplex 12.7.1
(with its default settings).  
FPCPA and CPBA may be run with a fixed proximity parameter, we use here $\mu=1.0$ which suits for well-scaled problems (FLA does not need the proximity parameter).
We ran FDSA with the rule \eqref{fdsa-mu}. 
As the sequence $\{\mu_k\}$ is decreasing, we consider
 a small positive constant $\mu_{\inf}$ and set  
$$
\mu_k=\max\lbrack\mu_{\inf}, \gamma_{k-1}\rbrack,\;k\ge 1,\;\mu_{\inf}=10^{-10}\|g^0\|.
$$ 
Our implementation of CPBA uses at each step $k$ a sequence of quadratic subproblemsfor $j=1,\ldots$
$$
z^{k,j}=\arg\min\limits_{(x,r)\in S\times\rit}\left\{\fmodel_{\bundleset_{k,j}}(x)+\frac{\mu}{2}\|x-\hat x^k\|^2\right\},
$$
where $\hat x^0=x^0$ and $\hat x^{k+1}=z^{k,j}$ if
$$
f(z^{k,j})\le f(\hat x^k)-\sigma\lbrack f(\hat x^k)-\fmodel_{\bundleset_{k,j}}(z^{k,j})\rbrack,
$$
in which case we have a descent step, otherwise a null step. In our experiments, we set $\sigma=0.5$.
Since they are nonsmooth unconstrained problems, we consider 
an input parameter $f_{\inf}^{(0)}$ to cope with the assumption of compactness of $S$. Hence, in Algorithm~\ref{fla} the lower bound is computed as 
$$\flow^k=\min\limits_{(x,r)\in\rit^{n+1}} \{\fmodel_{\bundleset_k}(x):\;f_{\inf}^{(0)}\le \fmodel_{\bundleset_k}(x)\},$$
and we add the constraint $f_{\inf}^{(0)}\le \fmodel_{\bundleset_k}(x)$ to all the quadratic subproblems to be consistent with our development. 
There are many tricks to avoid computing $\flow^k$ at each step, e.g. \cite{Frangioni2020,Kiwiel1995,OlSo2016} but for simplicity, it is updated as indicated above.
For all the test problems, we set $\kappa=0.8$ in \eqref{level-def} and $f_{\inf}^{(0)}=-10$ 
except for the problems 8 and 9 for which it takes the values $-100$ and $0$ respectively.
\begin{table}[]
\footnotesize
\caption{\small Test problems}
\begin{center}
\begin{tabular}{|c|c|c|c|}
\hline
 Problem & Name & $n$ & $f^*$  \\
\hline
1 & CB2 & 2 & 1.952224 \\
\hline
2 & CB3 & 2 & 2  \\
\hline
3 & DEM & 2 & -3  \\
\hline
4 & QL & 2 & 7.2 \\
\hline
5 & LQ & 2 & -$\sqrt{2}$ \\
\hline
6 & Mifflin1 & 2 & -1\\
\hline
7 & Mifflin2 & 2 & -1\\
\hline
8 & Rosen-Suzuki & 4 & -44\\
\hline
9 & Shor & 5 & 22.600162\\
\hline
10 & Maxquad & 10 & -0.841408\\
\hline
11 & Maxq & 20 & 0 \\
\hline
12  & Maxl & 20 & 0 \\
\hline
13  & Goffin & 50 & 0 \\
\hline
14  & MxHilb & 50 & 0 \\
\hline
15  & L1Hilb & 50 & 0 \\
\hline
\end{tabular}
\label{tpb}
\end{center}
\end{table}
The maximum number of steps allowed for all the algorithms (number of descent steps in CPBA) is set to $500$. 
With the given optimal functions values,
we stop the algorithms on the same basis, when
$$
\fbest^k-f^*\le 10^{-6}(1+|\fbest^k|).
$$
\begin{table}[]
\scriptsize
\caption{\small Computational results} 
\begin{center}
\begin{tabular}{|c||c|c|c||c|c||}
  \hline
 & \multicolumn{3}{|c||}{CPBA} & \multicolumn{2}{|c||}{FPCPA1}  \\  
\cline{2-6} \raisebox{1.5ex}[0cm][0cm]{Pb}  & $\#k$ & $\#fg$ & $f-f^*$ & $\# fg$ & $f-f^*$ \\
\hline
  1 & 8 & 22 & 9.35E-07 & 23 &  1.64E-06 \\
\hline
2 & 7 & 14 & 3.51E-07 & 12 & 4.54E-08 \\
\hline
3  & 4 & 7 & 3.08E-09 &  8 & 3.84E-09 \\
\hline
4  & 8 & 20 & 2.57E-06 & 27 &  1.97E-06 \\
\hline
5  & 4 & 8 & 1.29E-07 &  6 & 1.75E-06  \\
\hline
6  & 9 & 27 & 6.70E-07 & 20  & 4.50E-07 \\
\hline
7  & 8 & 22 & 1.13E-06 & 22 & 3.74E-07 \\
\hline
  8  & 9 & 40 & 3.90E-05 & 40 &  3.74E-05 \\
\hline
  9  & 10 & 43 & 1.89E-05 & 43 & 1.45E-05 \\
\hline
 10 & 14 & 127 & 3.94E-07 & 209 & 9.56E-07 \\
\hline
11  & 78 & 456 & 9.59E-07 & 269 & 7.60E-07  \\
\hline
12  & 210 & 231 & 6.36E-08 & 81 &  5.89E-09 \\
\hline
13  & 25 & 69 & 2.47E-10 & 87 &  8.36E-10 \\
\hline
14  & 500$^\dag$ & 504 & 1.50E-04 & 264 &  9.72E-07   \\
\hline
15  & 161 & 433 & 9.76E-07 & 62  & 8.40E-07 \\
  \hline
\end{tabular}
\\
\begin{tabular}{|c||c|c||c|c||}
\hline 
\multicolumn{5}{|c||}{} \\ \hline
 & \multicolumn{2}{|c||}{FLA1} & \multicolumn{2}{|c||}{FDSA1}  \\  
\cline{2-5} \raisebox{1.5ex}[0cm][0cm]{Pb}  & $\#fg$ & $f-f^*$ & $\# fg$ & $f-f^*$ \\
\hline
  1  & 18  & 7.46E-07 & 22 & 2.53E-06\\
\hline
2  & 16 & 4.52E-07   &13 & 3.13E-07\\
\hline
3  &  11  &  9.55E-07  & 11 & 2.58E-06\\
\hline
4   & 17  &  6.60E-06    & 19 & 1.47E-06\\
\hline
5  &  11 & 2.51E-07   & 7 & 3.79E-08\\
\hline
6  & 21  & 1.95E-06  & 16 & 1.39E-08\\
\hline
7  &  27 & 6.73E-07  & 17 & 1.91E-06\\
\hline
  8  & 70  &  2.47E-05  & 48 & 9.41E-06\\
\hline
  9  & 59  &  1.50E-05  & 41 & 2.25E-05\\
\hline
 10 & 204 & 1.34E-06  & 202 & 1.29E-06\\
\hline
11  & 231 & 9.11E-07  & 77 & 2.97E-07\\
\hline
12  &  48 & 5.12E-07  & 8 & 2.99E-09\\
\hline
13  &  59 &  1.86E-10 & 50 & 5.81E-12\\
\hline
14   & 19  &  2.71E-07  & 8 & 1.52E-07 \\
\hline
15  &  26  & 8.72E-07   & 8 & 5.17E-07\\
\hline
\end{tabular}
\label{results}
\end{center}
$\dag$ maximum number of $k$-steps (500) reached. 
\end{table}
We report on Table~\ref{results}, the number of steps (column $\#k$) for CPBA, the number of steps is the same as
the number of calls to $f$-oracle (column $\#fg$ which also indicates the number of steps of all the algorithms except CPBA)
and the absolute difference between $f$, the best function value found at stop and the optimal value $f^*$.  
These experiments show an improvement of the first two proposed algorithms over the classical proximal bundle algorithm since
both solve all the test problems within the maximum number of steps allowed to the contrary of the latter.
The rule \eqref{fdsa-mu} (which gives a decreasing sequence of proximity parameters) seems to be effective with Algorithm~\ref{fast-dbl} which compares favorably to 
the other algorithms on a majority of the test problems in terms of number of calls to the oracle. We expect further improvement from a more sophisticated management
of the proximity parameter in this algorithm.
\section{Concluding remarks}
We developed new algorithms for nonsmooth convex problems in the line of our previous approach in \cite{Ouorou2019} based of fast
gradient methods for smooth optimization. The limited experiments to get a first look at their performances is encouraging. Numerical experiments on large scale
problems are needed to confirm these performances including the benefit analysis of the momentum term by G\"uler. Another question we would
like to investigate is whether the use of non Euclidean entropy-like distances may be beneficial in the present setting as it is  
for the classical proximal bundle algorithms on certain convex problems. See 
the recent synthesis in \cite{Teboulle2018} exposing the benefits and limitations of the non Euclidean proximal framework.
\section*{Acknowledgements}
  I'm very thankful to Walid Ben Ameur for his comments on a previous version of the paper.

\bibliographystyle{cas-model2-names}
\bibliography{fast-prox-refs}

\begin{thebibliography}{16}
\expandafter\ifx\csname natexlab\endcsname\relax\def\natexlab#1{#1}\fi
\providecommand{\url}[1]{\texttt{#1}}
\providecommand{\href}[2]{#2}
\providecommand{\path}[1]{#1}
\providecommand{\DOIprefix}{doi:}
\providecommand{\ArXivprefix}{arXiv:}
\providecommand{\URLprefix}{URL: }
\providecommand{\Pubmedprefix}{pmid:}
\providecommand{\doi}[1]{\href{http://dx.doi.org/#1}{\path{#1}}}
\providecommand{\Pubmed}[1]{\href{pmid:#1}{\path{#1}}}
\providecommand{\bibinfo}[2]{#2}
\ifx\xfnm\relax \def\xfnm[#1]{\unskip,\space#1}\fi
\bibitem[{Attouch and Cabot(2018)}]{AtCa2018}
\bibinfo{author}{Attouch, H.}, \bibinfo{author}{Cabot, A.},
  \bibinfo{year}{2018}.
\newblock \bibinfo{title}{Convergence rates of inertial forward-backward
  algorithms}.
\newblock \bibinfo{journal}{SIAM J. Optim.} \bibinfo{volume}{28(1)},
  \bibinfo{pages}{849--874}.
\bibitem[{Correa and Lemar\'echal(1993)}]{CoLe1993}
\bibinfo{author}{Correa, R.}, \bibinfo{author}{Lemar\'echal, C.},
  \bibinfo{year}{1993}.
\newblock \bibinfo{title}{New variants of bundle methods}.
\newblock \bibinfo{journal}{Mathematical Programming} \bibinfo{volume}{62},
  \bibinfo{pages}{261--275}.
\bibitem[{Drori and Teboulle(2014)}]{DroriTeboulle2014}
\bibinfo{author}{Drori, Y.}, \bibinfo{author}{Teboulle, M.},
  \bibinfo{year}{2014}.
\newblock \bibinfo{title}{Performance of first-order methods for smooth convex
  minimization: a novel approach}.
\newblock \bibinfo{journal}{Mathematical Programming} \bibinfo{volume}{145},
  \bibinfo{pages}{451--482}.
\bibitem[{Frangioni(2020)}]{Frangioni2020}
\bibinfo{author}{Frangioni, A.}, \bibinfo{year}{2020}.
\newblock \bibinfo{title}{Standard bundle methods: Untrusted models and
  duality}, in: \bibinfo{editor}{Bagirov, A.}, \bibinfo{editor}{Gaudioso, M.},
  \bibinfo{editor}{Karmitsa, N.}, \bibinfo{editor}{M{\"a}kel{\"a}, M.} (Eds.),
  \bibinfo{booktitle}{Numerical nonsmooth optimization}.
  \bibinfo{publisher}{Springer}. chapter~\bibinfo{chapter}{3}, pp.
  \bibinfo{pages}{61--116}.
\bibitem[{G\"uler(1992)}]{Guler1992}
\bibinfo{author}{G\"uler, O.}, \bibinfo{year}{1992}.
\newblock \bibinfo{title}{New proximal point algorithm for convex
  minimization}.
\newblock \bibinfo{journal}{SIAM J. On Optimization} \bibinfo{volume}{2(4)},
  \bibinfo{pages}{649--664}.
\bibitem[{Hiriart-Urruty and Lemar\'echal(1993)}]{HUL1993}
\bibinfo{author}{Hiriart-Urruty, J.B.}, \bibinfo{author}{Lemar\'echal, C.},
  \bibinfo{year}{1993}.
\newblock \bibinfo{title}{Convex analysis and minimization algorithms}.
\newblock \bibinfo{howpublished}{Springer, Berlin}.
\bibitem[{Kim and Fessler(2016)}]{KimFessler2016}
\bibinfo{author}{Kim, D.}, \bibinfo{author}{Fessler, J.}, \bibinfo{year}{2016}.
\newblock \bibinfo{title}{Optimized first-order methods for smooth convex
  minimization}.
\newblock \bibinfo{journal}{Mathematical Programming} \bibinfo{volume}{159},
  \bibinfo{pages}{81--107}.
\bibitem[{Kim and Fessler(2018)}]{KimFessler2018}
\bibinfo{author}{Kim, D.}, \bibinfo{author}{Fessler, J.}, \bibinfo{year}{2018}.
\newblock \bibinfo{title}{Generalizing the optimized gradient method for smooth
  convex minimization}.
\newblock \bibinfo{journal}{SIAM J. On Optimization} \bibinfo{volume}{28(2)},
  \bibinfo{pages}{1920--1950}.
\bibitem[{Kiwiel(1995)}]{Kiwiel1995}
\bibinfo{author}{Kiwiel, K.}, \bibinfo{year}{1995}.
\newblock \bibinfo{title}{Proximal level bundle methods for convex
  nondifferentiable optimization, saddle-point problems and variational
  inequalities}.
\newblock \bibinfo{journal}{Mathematical Programming} \bibinfo{volume}{69},
  \bibinfo{pages}{89--109}.
\bibitem[{Lemar\'echal et~al.(1995)Lemar\'echal, Nemirovskii and
  Nesterov}]{LNN95}
\bibinfo{author}{Lemar\'echal, C.}, \bibinfo{author}{Nemirovskii, A.},
  \bibinfo{author}{Nesterov, Y.}, \bibinfo{year}{1995}.
\newblock \bibinfo{title}{New variants of bundle methods}.
\newblock \bibinfo{journal}{Mathematical Programming} \bibinfo{volume}{69},
  \bibinfo{pages}{111--147}.
\bibitem[{Luk\u{s}an(2000)}]{LukV2000}
\bibinfo{author}{Luk\u{s}an, L.~Vl\u{c}ek, J.}, \bibinfo{year}{2000}.
\newblock \bibinfo{title}{Test problems for nonsmooth unconstrained and
  linearly constrained optimization}.
\newblock \bibinfo{journal}{Technical report 798, Institute of Computer
  Science, Academiy of Sciences of the Czech Republic} .
\bibitem[{Nesterov(1983)}]{Nesterov1983}
\bibinfo{author}{Nesterov, Y.}, \bibinfo{year}{1983}.
\newblock \bibinfo{title}{A method for solving a convex programming problem
  with convergence rate $\text{O}(1/k^2)$}.
\newblock \bibinfo{journal}{Dokl. Akad. Nauk SSSR} \bibinfo{volume}{269},
  \bibinfo{pages}{543--547}.
\bibitem[{Nesterov(2004)}]{Nesterov2004}
\bibinfo{author}{Nesterov, Y.}, \bibinfo{year}{2004}.
\newblock \bibinfo{title}{Introductory lectures on convex programming. \text{A}
  basic course.}
\newblock \bibinfo{howpublished}{Kluwer Boston}.
\bibitem[{de~Oliveira and Solodov(2016)}]{OlSo2016}
\bibinfo{author}{de~Oliveira, W.}, \bibinfo{author}{Solodov, M.},
  \bibinfo{year}{2016}.
\newblock \bibinfo{title}{A doubly stabilized bundle method for nonsmooth
  convex optimization}.
\newblock \bibinfo{journal}{Mathematical Programming} \bibinfo{volume}{156(1)},
  \bibinfo{pages}{125--159}.
\bibitem[{Ouorou(2020)}]{Ouorou2019}
\bibinfo{author}{Ouorou, A.}, \bibinfo{year}{2020}.
\newblock \bibinfo{title}{Proximal bundle algorithms for nonsmooth convex
  optimization via fast gradient smooth methods}.
\newblock \bibinfo{journal}{arXiv preprint arxiv:2003.03437} .
\bibitem[{Teboulle(2018)}]{Teboulle2018}
\bibinfo{author}{Teboulle, M.}, \bibinfo{year}{2018}.
\newblock \bibinfo{title}{A simplified view of first order methods for
  optimization}.
\newblock \bibinfo{journal}{Mathematical Programming} \bibinfo{volume}{170},
  \bibinfo{pages}{67--96}.

\end{thebibliography}

\end{document}